\documentclass[a4paper,11pt,leqno]{article}

\usepackage{amsmath}
\usepackage{amsfonts}
\usepackage{amssymb}
\usepackage{latexsym}
\usepackage{epsfig}
\usepackage{graphicx}
\usepackage{oldgerm}
\usepackage{theorem}
\usepackage[mathscr]{eucal}
\usepackage{bm}

\def\Re{{\rm Re} \,}
\def\Im{{\rm Im} \,}
\def\arg{{\rm arg}\,}
\def\hcap{{\rm hcap}}
\def\supp{{\rm supp}}

\def\rO{{\rm O}}
\def\Var{{\rm Var}}

\def\R{\mathbb{R}}
\def\N{\mathbb{N}}
\def\C{\mathbb{C}}
\def\W{\mathbb{W}}
\def\H{\mathbb{H}}
\def\O{\mathbb{O}}
\def\E{\mathbb{E}}

\def\x{\bm{x}}

\def\bLambda{\bm{\Lambda}}
\def\X{\bm{X}}
\def\Y{\bm{Y}}
\def\1{{\bf 1}}
\def\U{\bm{U}}

\def\cF{{\CMcal F}}
\def\cM{{\CMcal M}}
\def\cS{{\CMcal S}}
\def\cC{{\CMcal C}}
\def\cH{{\CMcal H}}
\def\cE{{\CMcal E}}

\def\sH{{\sf H}}
\def\sU{{\sf U}}
\def\sD{{\sf D}}
\def\sS{{\sf S}}
\def\sF{{\sf F}}

\def\law={\stackrel{\rm (law)}{=}}

\def\P{\mathbb{P}}
\def\E{\mathbb{E}}

\def\bP{{\bf P}}

\def\Urho{\underline{\rho}}
\def\Us{\underline{s}}

\def\SLE{1\mathrm{SLE}}
\def\NSLE{N\mathrm{SLE}}
\def\GFF{\mathrm{GFF}}

\theorembodyfont{\itshape}

\newtheorem{thm}{Theorem}[section]
\newtheorem{lem}[thm]{Lemma}

\newtheorem{prop}[thm]{Proposition}

\newtheorem{df}[thm]{Definition}

\newcommand{\SSC}[1]{\section{#1}\setcounter{equation}{0}}
\newcommand{\qed}{\hbox{\rule[-2pt]{3pt}{6pt}}}

\begin{document}

\title{\bf 
Gaussian free fields coupled with multiple SLEs \\
driven by stochastic log-gases 
\footnote{
This manuscript was prepared for the proceeding of
the workshop of
the 12th Mathematical Society of Japan, Seasonal Institute (MSJ-SI), 
`Stochastic Analysis,
Random Fields and Integrable Probability', held 
at Kyushu University, Fukuoka, Japan, 
July 31--August 9, 2019. 
The proceeding will be published in
\textit{Advanced Studies in Pure Mathematics 87, 2021, 
Stochastic Analysis, Random Fields and 
Integrable Probability -- Fukuoka 2019}.
}
}
\author{
Makoto Katori
\footnote{
Department of Physics,
Faculty of Science and Engineering,
Chuo University, 
Kasuga, Bunkyo-ku, Tokyo 112-8551, Japan;
e-mail: katori@phys.chuo-u.ac.jp} 
\, and \,
Shinji Koshida
\footnote{
Department of Physics,
Faculty of Science and Engineering,
Chuo University, 
Kasuga, Bunkyo-ku, Tokyo 112-8551, Japan;
e-mail: koshida@phys.chuo-u.ac.jp} 
}
\date{14 July 2021}
\pagestyle{plain}
\maketitle

\begin{abstract}
Miller and Sheffield introduced the notion of
an imaginary surface as an equivalence class
of pairs of simply connected proper subdomains of
$\mathbb{C}$ and Gaussian free fields (GFFs)
on them under the conformal equivalence.
They considered the situation in which
the conformal maps are given by
a chordal Schramm--Loewner evolution (SLE). 
In the present paper, we construct 
GFF-valued processes on $\mathbb{H}$ (the upper half-plane) 
and $\mathbb{O}$ (the first orthant of $\mathbb{C}$)
by coupling a GFF with a multiple SLE evolving in time on each domain.
We prove that a GFF on $\H$ and $\mathbb{O}$ is locally coupled 
with a multiple SLE if the multiple SLE is driven 
by the stochastic log-gas called the Dyson model defined 
on $\mathbb{R}$ and the Bru--Wishart process 
defined on $\mathbb{R}_+$, respectively.
We obtain pairs of time-evolutionary domains 
and GFF-valued processes.


Schramm--Loewner evolution,
Multiple SLE, 
Stochastic log-gases, 
Dyson model, 
Bru--Wishart process

\end{abstract}

\SSC{Introduction} \label{sec:introduction}

The present study is motivated by the recent work
by Sheffield on the quantum gravity zipper and the AC geometry \cite{She16}
and a series of papers by Miller and Sheffield on the 
imaginary geometry \cite{MS16a,MS16b,MS16c,MS17}.
In both of them, a {\it Gaussian free field} (GFF) on 
a simply connected proper subdomain $D$ of the complex plane $\C$ 
(see, for instance, \cite{She07}) is coupled
with a {\it Schramm--Loewner evolution} (SLE) 
\cite{Sch00,LSW04,Law05} driven
by a Brownian motion moving 
on the boundary $\partial D$, 
or its variant called an SLE($\kappa, \Urho$).

Consider a simply connected domain $D \subsetneq \C$
and write $\cC_{\rm c}^{\infty}(D)$ for the space of 
real smooth functions
on $D$ with compact support.
Assume $h \in \cC_{\rm c}^{\infty}(D)$ and 
consider a smooth vector field
$e^{\sqrt{-1} (h/\chi + \theta)}$ with parameters $\chi, \theta \in \R$.
Then, the {\it flow line} along this vector field, 
$\eta : (0, \infty) \ni t \mapsto \eta(t) \in D$,
starting from 
$\lim_{t \to 0} \eta(t) =: \eta(0)=x \in \partial D$ 
is defined (if exists) as the 
solution of the ordinary differential equation (ODE) 
\begin{equation}
\frac{d \eta(t)}{dt} = e^{\sqrt{-1}\{h(\eta(t))/\chi + \theta\}},
\quad t \geq 0, \quad \eta(0) = x.
\label{eqn:flow1}
\end{equation}
Let $\widetilde{D} \subsetneq \C$ be another 
simply connected domain and
consider a conformal map $\varphi: \widetilde{D} \to D$.
Then, we define the pull-back of the flow line $\eta$ by $\varphi$ as
$\widetilde{\eta}(t)=(\varphi^{-1} \circ \eta)(t)$.
That is, 
$\varphi(\widetilde{\eta}(t))=\eta(t)$, 
and the derivatives with respect to $t$ of both sides of this equation 
gives 
$\varphi'(\widetilde{\eta}(t)) d \widetilde{\eta}(t)/dt=d \eta(t)/dt$
with $\varphi'(z) := d \varphi(z)/dz$.
We use the polar coordinates
$\varphi'(\cdot)=|\varphi'(\cdot)| e^{\sqrt{-1} \arg \varphi'(\cdot)}$,
where $\arg \zeta$ of $\zeta \in \C$ is a priori 
defined up to additive multiples of $2 \pi$, and hence,
we have 
$d \widetilde{\eta}(t)/dt=
e^{\sqrt{-1}\{ (h \circ \varphi 
- \chi \arg \varphi')(\widetilde{\eta}(t))/\chi + \theta\}}
/|\varphi'(\widetilde{\eta}(t))|, t \geq 0$.
If we perform a time change $t \to \tau=\tau(t)$ by putting
$t=\int_0^{\tau} ds/|\varphi'(\widetilde{\eta}(s))|$ and 
$\widehat{\eta}(t):=\widetilde{\eta}(\tau(t))$, then 
the above equation becomes
\[
\frac{d \widehat{\eta}(t)}{dt}
= e^{\sqrt{-1}\{ (h \circ \varphi 
- \chi \arg \varphi')(\widehat{\eta}(t))/\chi + \theta\}},
\quad t \geq 0.
\]
Since a time change preserves the image of a flow line,
we can identify $h$ on $D$ and 
$h \circ \varphi - \chi \arg \varphi'$ 
on $\widetilde{D}=\varphi^{-1}(D)$.
In \cite{She16,MS16a,MS16b,MS16c,MS17}, such a flow line
is considered also in the case that $h$ is given by 
an instance of a GFF defined as follows.
\begin{df}
\label{thm:def_GFF}
Let $D\subsetneq \mathbb{C}$ be a simply connected 
domain and $H$ be 
a GFF on $D$ with zero boundary condition
(constructed in Section \ref{sec:GFF}). 
A GFF on $D$ is a random distribution $h$ 
of the form $h=H+u$, where $u$ is a 
deterministic harmonic function on $D$.
\end{df} 
\noindent
Since a GFF is not function-valued, 
but it is a {\it distribution-valued random field}
(see Remark 4.1 in Section \ref{sec:GFF}),
the ODE in the form (\ref{eqn:flow1}) 
no longer makes sense mathematically 
in the classical sense. 
Using the theory of SLE, however, 
the notion of flow lines was generalized as follows.

Consider the collection
\[
\sS:=\left\{(D,h) \Bigg| 
\substack{\, \ D\subsetneq \mathbb{C}:\ 
\mbox{\small simply connected} \\ 
h:\ \mbox{\small GFF on $D$}} \right\}.
\]
Fixing a parameter $\chi \in \R$, 
we define the following equivalence relation in $\sS$.
\begin{df}
\label{thm:def_IS}
Two pairs $(D,h)$ and
$(\widetilde{D}, \widetilde{h}) \in \sS$ are 
{\it equivalent} if there exists a conformal map 
$\varphi: \widetilde{D} \to D$ and
$\widetilde{h} 
\law= h \circ \varphi - \chi \arg \varphi^{\prime}$. 
In this case, we write 
$(D, h) \sim (\widetilde{D}, \widetilde{h})$.
\end{df}
\noindent
We call each element belonging to $\sS/\sim$ 
an {\it imaginary surface} \cite{MS16a} 
(or an {\it AC surface} \cite{She16}). 
That is, in this equivalence class, 
a conformal map $\varphi$ causes not only 
a coordinate change of a GFF as $h \mapsto h \circ \varphi$
associated with changing the domain of definition of the field as
$D \mapsto \varphi^{-1}(D)$,
but also an addition of a deterministic 
harmonic function $-\chi \arg \varphi'$ to the field.
Notice that this definition depends on
one parameter $\chi \in \R$.

As will be explained in Section \ref{sec:GFF}, 
each instance $H$ of a GFF with zero boundary condition
depends on the choice of 
a complete orthonormal system (CONS)
of a Hilbert space
starting from which a GFF is constructed.
The probability law of a zero-boundary GFF 
is, however, independent 
of such construction and uniquely determined.

Consider the case in which $D$ is the
upper half-plane $\H :=\{z \in \C: \Im z > 0 \}$
with $\partial \H=\R \cup \{\infty \}$.
Let $(B(t))_{t \geq 0}$ be a one-dimensional
standard Brownian motion starting from the origin 
defined on a probability space 
$(\Omega^{\SLE}, \cF^{\SLE}, \P^{\SLE})$ and
adapted to a filtration $(\cF^{\SLE}_{t})_{t\geq 0}$.
We consider the chordal SLE($\kappa$) 
driven by $(\sqrt{\kappa} B(t))_{t \geq 0}$ on
$S :=\R$ with $\kappa>0$ \cite{Sch00,LSW04,Law05}, 
associated to which we obtain a random curve 
(called a {\it chordal SLE($\kappa$) curve}) 
parameterized by time, 
$\eta: (0, \infty) \ni t \mapsto \eta(t) \in \H$, 
such that $\lim_{t \to 0} \eta(t) =:\eta(0)=0$, 
$\lim_{t \to \infty} \eta(t)=\infty$. 
At each time $t>0$, 
let $\eta(0, t] :=\{\eta(s) : s \in (0, t]\}$ and
we write $\H^{\eta}_{t}$ for 
the unbounded component of $\H\setminus\eta (0,t]$.
Then, the chordal SLE($\kappa$) gives a conformal map from 
$\H^{\eta}_t$ to $\H$. It is also known that, 
if $\kappa\in (0,4]$, then $\eta (0,t]$ is almost surely 
a simple curve at each $t>0$ and, hence, 
$\H^{\eta}_{t}=\H\setminus \eta(0,t]$. 
In this paper, we will write the chordal SLE($\kappa$) as 
$(g_{\H^{\eta}_t})_{t \geq 0}$.
Let $H(\cdot)$ be a GFF on $\H$ 
with zero boundary condition on $\R$ 
that is defined on a probability space 
$(\Omega^{\GFF},\cF^{\GFF},\P^{\GFF})$. 
To couple the SLE and the GFF, 
we introduce a probability space 
$(\Omega, \cF, \P)=(\Omega^{\GFF}\times \Omega^{\SLE}, 
\cF^{\GFF} \vee \cF^{\SLE}, \P^{\GFF} \otimes \P^{\SLE})$
and extend the SLE and the GFF onto this probability space. 
Then, the SLE is adapted to the filtration 
$(\cF_{t})_{t\geq 0}$ defined by 
$\cF_{t}=\{\emptyset, \Omega^{\GFF} \} \vee \cF^{\SLE}_{t}$.
Instead of $H(\cdot)$ itself, we consider the following GFF on $\H$
by adding a deterministic harmonic function,
\begin{equation}
h(\cdot) := H(\cdot) - \frac{2}{\sqrt{\kappa}} \arg (\cdot).
\label{eqn:arg1}
\end{equation}
Notice that $\arg(\cdot) = \Im \log (\cdot)$
and the real and imaginary parts of 
a complex analytic function are harmonic.
Hence, the random distribution (\ref{eqn:arg1}) 
is in fact a GFF in the sense of Definition \ref{thm:def_GFF}.
Given $\kappa >0$ for the SLE($\kappa$),
we fix the parameter $\chi$ as
$\chi = 2/\sqrt{\kappa}-\sqrt{\kappa}/2$. 
Note that the well-known relation between $\kappa$
and the {\it central charge} $c$ of 
conformal field theory 
is simply expressed using 
the present parameter $\chi$ as $c=1-6 \chi^2$
(see, for instance, \cite[Eq.(6)]{BB02}).
Let
\[
f_{\H^{\eta}_t} := g_{\H^{\eta}_t} - \sqrt{\kappa} B(t)
=\sigma_{-\sqrt{\kappa} B(t)} \circ g_{\H^{\eta}_t},
\]
where $\sigma_s$ denotes
the translation by $s \in \R$; 
$\sigma_s(z)=z+s$, $z \in \H$. 
Let $A\subset \H$ be an open set 
and take an $(\cF_{t})_{t\geq 0}$-stopping time
\begin{equation*}
\tau_A := \inf \left\{t \geq 0 \ \Big|
\eta (0, t] \cap A \not= \emptyset \right\}.
\label{eqn:tauA}
\end{equation*}
Let $\tau$ be any $(\cF_{t})_{t \geq 0}$-stopping time
such that $\tau \leq \tau_A$ a.s.
Then, we can prove the following equality in probability
\cite[Theorem 1.1, Lemma 3.11]{MS16a} 
(see also \cite[Lemma 6.1]{Dub09}); 
for any $f\in \cC_{\rm c}^{\infty}(\H)$ such that $\supp (f)\subset A$,
\begin{equation}
(h, f)
\law= (h \circ f_{\H^{\eta}_{\tau}} 
- \chi \arg f_{\H^{\eta}_{\tau}}^{\prime}, f) \quad
\mbox{under $\P$},  
\label{eqn:chi_equivalence1}
\end{equation}
where the pairing $(\cdot, \cdot)$ is defined 
by (\ref{eqn:inner_product}) below. 
We comment that, 
due to the conformal invariance of 
a zero-boundary GFF
(see Section \ref{sec:conformal_invariance} below), 
for an instance of the SLE($\kappa$), the random distribution 
$h \circ f_{\H^{\eta}_{t}} - \chi \arg f_{\H^{\eta}_{t}}^{\prime}$
is a GFF on $\H^{\eta}_{t}$ in the sense of Definition 
\ref{thm:def_GFF}. 
Notice that pairs $(\H, h)$
and $(\H^{\eta}_t, h \circ f_{\H^{\eta}_t} - \chi \arg f_{\H^{\eta}_t}')$
with (\ref{eqn:arg1}) are equivalent 
in the sense of Definition \ref{thm:def_IS}.
In other words, an imaginary surface 
whose representative is given by $(\H, h)$
is constructed as
a pair of a time-evolutionary domain, 
$f^{-1}_{\H^{\eta}_t}(\H)=\H^{\eta}_t$, 
$t \geq 0$,
and a GFF-valued process, 
$h \circ f_{\H^{\eta}_t}-\chi \arg f_{\H^{\eta}_t}', t \geq 0$
defined on it.
With the establishment of the equality (\ref{eqn:chi_equivalence1})
we say that the \textit{local coupling 
between a GFF and an SLE} 
is constructed (see \cite{Dub09,She16,MS16a}
for lifting the local coupling to the `global' one). 
It was proved \cite{Dub09,MS16a} that, 
under the coupling between a GFF and an SLE, 
the SLE-curve is a deterministic functional of the GFF. 
By virtue of it, in \cite{MS16a}, 
the authors referred to an SLE($\kappa$) curve 
as a flow line of the GFF $h$.

Here, first we consider 
the case in which the conformal maps
are generated by a multiple Loewner equation associated with
a multi-slit.
Let $N \in \N :=\{1,2, \dots\}$ and suppose
that we have $N$ slits 
$\eta_i =\{\eta_i(t): t \in (0, \infty)\} \subset \H$, 
$1 \leq i \leq N$,
which are simple curves, disjoint with each other, 
$\eta_i \cap \eta_j = \emptyset, i \not= j$,
starting from $N$ distinct points 
$\lim_{t \to 0} \eta_i(t) =: \eta_i(0)$ on $\R$;
$\eta_1(0) < \cdots < \eta_N(0)$, 
and all going to infinity; $\lim_{t \to \infty} \eta_i(t)=\infty$, 
$1 \leq i \leq N$. 
A {\it multi-slit} is defined as the union of them, 
$\bigcup_{i=1}^N \eta_i$, and 
\begin{align*}
\H^{\eta}_t &:= \H \setminus \bigcup_{i=1}^N \eta_i(0, t]
\ \
\mbox{for each $t > 0$ with $\H^{\eta}_0 := \H$}. 
\end{align*}
We write the time evolution of the conformal map
which transforms $\H^{\eta}_t$
to $\H$ at each time $t \geq 0$ 
under the hydrodynamic normalization
as $(g_{\H^{\eta}_t})_{t \geq 0}$
and call it a {\it multiple SLE}.
The images of the tips of the multi-slit 
$g_{\H^{\eta}_t}(\eta_i(t))$,
$1 \leq i \leq N$ exist as points on $\R$ for $t \geq 0$
and if we put $X^{\R}_i(t) := g_{\H^{\eta}_t}(\eta_i(t))$,
the multiple SLE $(g_{\H^{\eta}_t})_{t \geq 0}$ is given as
a unique solution 
of the following equation,
\begin{align}
\frac{d g_{\H^{\eta}_t}(z)}{dt}
&=\sum_{i=1}^N \frac{2}{g_{\H^{\eta}_t}(z)-X^{\R}_i(t)}, \quad t \geq 0,
\label{eqn:mSLE0}
\\
g_{\H^{\eta}_0}(z) &=z \in \H,
\nonumber
\end{align}
under a proper parameterizatoin of the multi-slit.
Here 
\[
\X^{\R}(t)=(X^{\R}_1(t), \dots, X^{\R}_N(t)) \in \R^N, 
\quad t \geq 0
\]
is called the {\it driving process} of the multiple SLE.

In the sequel, we will consider the case when 
$(\X^{\R}(t))_{t\geq 0}$ is a stochastic process 
defined on a probability space 
$(\Omega^{\NSLE}, \cF^{\NSLE}, \P^{\NSLE})$ and 
adapted to a filtration $(\cF^{\NSLE}_{t})_{t\geq 0}$. 
In this case, although it is not ensured that the solution 
generates a multi-slit depending on $\kappa$,
we can still find a family of domains 
$\H^{\eta}_{t}\subseteq\H$, $t \geq 0$
so that $g_{\H^{\eta}_{t}}: \H^{\eta}_{t} \to \H$ 
is a conformal map at each $t \geq 0$.
(See Remark 1.4 (1) in \cite{KK21}.)

We again consider a zero-boundary GFF $H$ on $\H$ defined 
on a probability space 
$(\Omega^{\GFF}, \cF^{\GFF}, \P^{\GFF})$ and introduce 
a coupled probability space 
\[
(\Omega, \cF, \P)=(\Omega^{\GFF} \times \Omega^{\NSLE},
\cF^{\GFF} \vee \cF^{\NSLE}, \P^{\GFF} \otimes \P^{\NSLE}). 
\]
Then, the multiple SLE and the GFF are naturally extended 
to $(\Omega, \cF, \P)$, and the multiple SLE is adapted to 
the filtration $(\cF_{t})_{t\geq 0}$ defined by 
$\cF_{t}=\{\emptyset, \Omega^{\GFF} \} \vee \cF^{\NSLE}_{t}$.

Regarding (\ref{eqn:arg1}) and (\ref{eqn:chi_equivalence1}), we see that 
$h \circ f_{\H^{\eta}_t}(\cdot)-\chi \arg f_{\H^{\eta}_t}'(\cdot)$ is
equal to 
\begin{align*}
& (H \circ \sigma_{-\sqrt{\kappa} B(t)} ) \circ g_{\H^{\eta}_t}(\cdot)
-\frac{2}{\sqrt{\kappa}} \arg(g_{\H^{\eta}_t}(\cdot)-\sqrt{\kappa} B(t))
-\chi \arg g_{\H^{\eta}_t}'(\cdot)
\nonumber\\
&
\quad
=
H \circ g_{\H^{\eta}_t}(\cdot)
-\frac{2}{\sqrt{\kappa}} 
\arg(g_{\H^{\eta}_t}(\cdot)-g_{\H^{\eta}_t}(\eta(t)))
-\chi \arg g_{\H^{\eta}_t}'(\cdot) 
\nonumber\\
& \hskip 7cm
\quad \mbox{under $\P$}, 
\end{align*}
$t \geq 0$, 
where the translation invariance of $H$ was used.
Motivated by this observation, 
we study the GFF-valued process defined by
\begin{align}
H_{\H}(\cdot, t)
:=& H \circ g_{\H^{\eta}_t}(\cdot) \\
& \,
-\frac{2}{\sqrt{\kappa}} \sum_{i=1}^N
\arg( g_{\H^{\eta}_t}(\cdot)-g_{\H^{\eta}_t}(\eta_i(t)))
-\chi \arg g_{\H^{\eta}_t}'(\cdot) 
\nonumber\\
=& H \circ g_{\H^{\eta}_t}(\cdot) 
\nonumber\\
& \,
-\frac{2}{\sqrt{\kappa}} \sum_{i=1}^N
\arg( g_{\H^{\eta}_t}(\cdot)-X^{\R}_i(t))
-\chi \arg g_{\H^{\eta}_t}'(\cdot) 
\nonumber
\label{eqn:GFF_mSLE1}
\end{align}
on $\H^{\eta}_t, t \geq 0$.
This process starts from 
\[
H_{\H}(\cdot,0)=H(\cdot)-
\frac{2}{\sqrt{\kappa}} \sum_{i=1}^{N} \arg (\cdot-x^{\R}_{i}),
\]
where we assume that $x^{\R}_1 < \cdots < x^{\R}_N$.
We let the boundary points evolve according to 
the stochastic process $(\X^{\R}(t))_{t \geq 0}$
starting from $\x^{\R}:=(x^{\R}_i)_{i=1}^N$.
At each time $t>0$, 
we consider the GFF $H+u_{t}$ on $\mathbb{H}$ where
$u_t(\cdot)=-(2/\sqrt{\kappa}) 
\sum_{i=1}^{N} \arg (\cdot- X^{\R}_{i}(t))$.
Then, the GFF $H_{\H}(\cdot,t)$ on $\H^{\eta}_{t}$ 
is defined by the property that 
$(\H,H+u_{t}) \sim (\H^{\eta}_{t}, H_{\H}(\cdot,t))$
in the sense of Definition \ref{thm:def_IS}.

A part of the main theorem in this paper (Theorem \ref{thm:main1}) 
is stated as follows.

\begin{thm}
\label{thm:main0}
Let $A\subset\H$ be an open subset and 
take an $(\cF_{t})_{t\geq 0}$-stopping time
\begin{equation}
\label{eqn:tau_A_multiple}
\tau_A := \inf \left\{t \geq 0 \ \Big| A \not\subseteq \H^{\eta}_{t} \right\}.
\end{equation}
Let $\tau$ be any $(\cF_{t})_{t\geq 0}$-stopping time 
such that $\tau \leq \tau_A$ a.s.
Then, for any $f\in \cC_{\rm c}^{\infty}(\H)$ 
such that $\supp (f)\subset A$,
\begin{equation}
(H_{\H}(\cdot, 0), f)
\law= (H_{\H}(\cdot, \tau), f )
\ \
\mbox{under $\P$},
\label{eqn:stationary0}
\end{equation}
if the driving process $(\X^{\R}(t))_{t \geq 0}$
is equal to the time changed version 
$\Y^{\R}(t)=(Y^{\R}_1(t), \dots, Y^{\R}_N(t))$, $t \geq 0$
of the Dyson model on $\R$
which solves the following system of
stochastic differential equations (SDEs)
with $\kappa >0$,
\begin{equation}
dY^{\R}_i(t)=\sqrt{\kappa} dB_i(t)
+ 4 \sum_{1 \leq j \leq N, j \not=i}
\frac{dt}{Y^{\R}_i(t)-Y^{\R}_j(t)},
\label{eqn:Dyson0}
\end{equation}
$t \geq 0, 1 \leq i \leq N$,
where 
$(B_i(t))_{t \geq 0}$
are mutually independent 
one-dimensional standard
Brownian motions
starting from
$B_i(0)= Y^{\R}_i(0) =: y^{\R}_i
=\eta_i(0), 1 \leq i \leq N$,
satisfying $y^{\R}_1 < \cdots < y^{\R}_N$.
\end{thm}

The Dyson model \cite{Dys62} is one of the most studied
stochastic log-gases in one dimension,
which is a dynamical version of the
one-parameter ($\beta=8/\kappa$) extension 
of the {\it Gaussian unitary ensemble (GUE)} 
of point processes 
studied in random matrix theory
\cite{For10,Kat15_Springer}.
It is also known that the multiple SLE driven by
$(\Y^{\R}(t))_{t\geq 0}$ is an example of multiple SLEs
that is defined in terms of an SLE partition function 
\cite{BBK05,Gra07}. 

It is possible that, for a certain choice of the open set 
$A\subset\H$ in Theorem \ref{thm:main0}, 
the stopping time $\tau$ has to be $\tau=0$.
Let us see that we can take $\tau>0$ a.s. 
for a generic choice of $A$.
For each $1 \leq i \leq N$, let $U_{i}$ be 
a neighborhood of $y^{\R}_{i}$ in $\overline{\H}$ and 
suppose that $U_{i}\cap U_{j}=\emptyset$ if $i\neq j$. 
We call such $U_{i}$, $1\leq i\leq N$
\textit{localization neighborhoods}. 
Then, we define the first exit time of the SLE 
from the union of localization neighborhoods, 
$\U :=\bigcup_{i=1}^N U_i$ by
\[
\tau_{\U}=\inf \left\{t \geq 0 \ \Big| 
(\H \setminus \U) 
\not\subseteq \H^{\eta}_{t} \right\}.
\]
Note that we have $\tau_{\bm{U}}>0$ a.s. 
since, otherwise, it contradicts the fact that 
$(\Y^{\R}(t))_{t\geq 0}$ is continuous in $t\geq 0$ \cite{CL97,GM14}.
Notice that we can take the stopping time $\tau$ in 
Theorem \ref{thm:main0} in such a way that 
$\tau \ge \tau_{\bm{U}}$ for some 
localization neighborhoods that are disjoint from 
$A$, if the closure of $A$ in $\overline{\H}$ 
does not contain any of $y_{i}^{\R}$, $1 \leq i \leq N$.
Therefore, we can take it so that $\tau>0$ a.s. 
in such a generic case.

As we have already pointed out, it is not clear if a solution 
of the multiple Loewner equation generates a multi-slit.
In our subsequent paper \cite{KK21}, 
we will prove under the assumption $\kappa\in (0,8]$ that, 
if the driving process is the time changed version 
$\Y^{\R}(t)=(Y^{\R}_1(t), \dots, Y^{\R}_N(t))$, $t \geq 0$
of the Dyson model as in Theorem \ref{thm:main0}, 
the solution generates a set of curves. 
We will also prove that the resulting curves $\{\eta_i\}_{i=1}^N$
are simple disjoint curves 
if $\kappa\in (0,4]$, self-intersecting if $\kappa\in (4,8)$, 
and space-filling if $\kappa=8$.

The process $(H_{\H}(\cdot,t))_{t\ge 0}$ 
is a generalization of 
$(h \circ f_{\H^{\eta}_t}-\chi \arg f_{\H^{\eta}_t}')_{t \geq 0}$
considered by Miller and Sheffield \cite{She16,MS16a}
as explained above.
The equality (\ref{eqn:chi_equivalence1}) 
\cite{Dub09,She16,MS16a} has been 
extended to the equality (\ref{eqn:stationary0}) 
in Theorem \ref{thm:main0}, 
which we think of as the local coupling 
between a GFF and a multiple SLE.

We will also construct another 
GFF-valued process in the first orthant in $\C$; 
$\O:=\{z \in \C : \Re z >0, \Im z >0 \}$. 
There, a GFF on $\O$, denoted as $H_{\O}(\cdot,0)$
is locally coupled with a multiple version of 
the {\it quadrant SLE} \cite{Tak14} defined on $\O$,
which is driven by a stochastic log-gas
defined on $S=\R_+ := \{x \in \R: x \geq 0 \}$.
This driving process is a dynamical version
of the one-parameter ($\beta=8/\kappa$) 
extension of the {\it chiral GUE} of point processes
with parameter $\nu \in [0, \infty)$ 
studied in random matrix theory \cite{KT04,For10}, and 
we call it the {\it Bru--Wishart process} 
in this paper \cite{Wis28,Bru91}.
We note that 
$(\H, H_{\H}(\cdot, 0)) \sim (\O, H_{\O}(\cdot, 0))$ 
in the sense of Definition \ref{thm:def_IS}.

Construction of such GFF-valued processes will be
meaningful for the study of multiple SLEs.
The main problem in defining a multiple SLE correctly
in $D \subsetneq \C$ may be how to find
a correct principle 
to choose a driving process $(\X^S(t))_{t \geq 0}$
defined on a part of the boundary $S \subset \partial D$
(e.g., conformal invariance, statistical mechanics consideration,
reparameterization invariance,
absolute continuity to the SLE with a single slit,
commutation relations) 
\cite{Car03,BBK05,KL07,Gra07,Dub07}.
In the present paper, we simply assume the form of
SDEs for $(\X^S(t))_{t \geq 0}$ as
\begin{equation}
d X^S_i(t)=\sqrt{\kappa} dB_i(t)+
F^S_i(\X^S(t)) dt, \quad t \geq 0,
\quad 1 \leq i \leq N,
\label{eqn:SDE0}
\end{equation}
where $(B_i(t))_{t \geq 0}, 1 \leq i \leq N$ are
mutually independent one-dimensional standard Brownian motions,
$\kappa > 0$, and
$F^S_i(\x) \in \cC^{\infty}
(S^N \backslash \bigcup_{j \not=k} \{x_j=x_k\})$,
$1 \leq i \leq N$, 
which  do not explicitly depend on $t$. 
Then, the equality (\ref{eqn:stationary0}) 
for a GFF-valued process on $D=\H$ 
determines the driving process $(\X^{\R}(t))_{t \geq 0}$ as
(a time change of) the Dyson model $(\Y^{\R}(t))_{t \geq 0}$.
That is, the local coupling between a GFF and a multiple SLE
provides a new scheme to choose a driving process for
a multiple SLE.

Notice again
that $\arg z$ in (\ref{eqn:arg1}) is the imaginary part of 
the complex analytic function $\log z$.
Sheffield studied another type of distribution-valued 
random field on $\H$ given by \cite{She16} 
\[
\widetilde{h}(\cdot) 
:= \widetilde{H}(\cdot)+ \frac{2}{\sqrt{\kappa}}
\Re \log (\cdot)
= \widetilde{H}(\cdot)+ \frac{2}{\sqrt{\kappa}} \log |\cdot|,
\]
where $\widetilde{H}(\cdot)$ is a free boundary GFF on $\H$
and found that $\widetilde{h}(\cdot)$ is 
coupled with a backward SLE in the context of 
quantum gravity \cite{DS11}.
This coupling was later generalized in \cite{KK20,Kos19} 
to the situations where backward multiple SLEs driven 
by stochastic log-gases play analogous roles as multiple SLEs 
did in the present work.

The present paper is organized as follows. 
We give brief reviews of stochastic log-gases
in one dimension in Section \ref{sec:log_gases}
and the SLE both for a single-slit and a multi-slit
in Section \ref{sec:multiple_SLE}. 
In Section \ref{sec:GFF}, we define
a GFF with zero boundary condition
on $D \subsetneq \C$ based on 
the Bochner--Minlos theorem.
The construction of GFF-valued processes
by locally coupling GFFs with multiple SLEs 
driven by specified stochastic log-gases on $S$
are given in Section \ref{sec:GFF_loggas}
for $(D, S)=(\H, \R)$ and $(\O, \R_+)$. 

\SSC{One-dimensional Stochastic Log-Gases} \label{sec:log_gases}
\subsection{Eigenvalue and singular-value processes}
\label{sec:matrix}

For $N \in \N$, 
let $\sH_N$ and $\sU_N$ be the space of 
$N \times N$ Hermitian matrices and
the group of $N \times N$ unitary matrices, respectively.
Consider complex-valued processes
$(M_{ij}(t))_{t \geq 0}, 1\leq i, j \leq N$ 
with the condition $\overline{M_{ji}(t)}=M_{ij}(t)$,
where $\overline{z}$ denotes the
complex conjugate of $z \in \C$.
We consider an $\sH_N$-valued process by
$M(t)=(M_{ij}(t))_{1 \leq i, j \leq N}$. 
For $S=\R$ and $\R_+$, define
the Weyl chambers as
$\W_N(S) := \{ \x=(x_1, \dots, x_N) \in S^N: x_1 < \cdots < x_N\}$,
and write their closures as
$\overline{\W_N(S)} = \{ \x \in S^N :
x_1 \leq \cdots \leq x_N\}$. 
For each $t \geq 0$, there exists 
$U(t)=(U_{ij}(t))_{1 \leq i, j \leq N} \in \sU_N$ such that
it diagonalizes $M(t)$ as
$U^{\dagger}(t) M(t) U(t) = {\rm diag}(\Lambda_1(t), \dots, 
\Lambda_N(t))$
with the eigenvalues $\{\Lambda_i(t)\}_{i=1}^N$ of $M(t)$,
where $U^{\dagger}(t)$ is the Hermitian conjugate of 
$U(t)$; $U^{\dagger}_{ij}(t)=\overline{U_{ji}(t)}, 1 \leq i, j \leq N$, 
and we assume
$\bLambda(t) :=(\Lambda_1(t), \dots, \Lambda_N(t))
\in \overline{\W_N(\R)}$, $t\geq 0$. 
For $d M(t):=(d M_{ij}(t))_{1 \leq i, j \leq N}$, define
a set of quadratic variations, 
\[
\Gamma_{ij, k \ell}(t):= \Big\langle (U^{\dagger} d M  U)_{ij}, 
(U^{\dagger} d M  U)_{k \ell} \Big\rangle_t,
\quad 1 \leq i, j, k, \ell \leq N, \quad t \geq 0.
\]
We write $\1_E$ for the indicator function
of an event $E$;
$\1_E=1$ if $E$ occurs,
and $\1_E=0$ otherwise.
The following is proved \cite{Bru89,KT04,Kat15_Springer}.
See Section 4.3 of \cite{AGZ10} for details of proof. 

\begin{prop}
\label{thm:Bru}
Assume that $(M_{ij}(t))_{t \geq 0}, 1 \leq i, j \leq N$ are
continuous semi-martingales.
The eigenvalue process $(\bLambda(t))_{t \geq 0}$ satisfies
the following system of SDEs,
\[
d \Lambda_i(t)=d \cM_i(t)+ d J_i(t), \quad t \geq 0, 
\quad 1 \leq i \leq N, 
\]
where $(\cM_i(t))_{t \geq 0}, 1 \leq i \leq N$
are martingales with quadratic variations
$\langle \cM_i, \cM_j \rangle_t
= \int_0^t \Gamma_{ii, jj}(s) ds$, 
and $(J_i(t))_{t \geq 0}, 1 \leq i \leq N$ are the
processes with finite variations given by
\[
dJ_i(t)= \sum_{j=1}^N 
\frac{\1_{\Lambda_i(t) \not= \Lambda_j(t)}}
{\Lambda_i(t)-\Lambda_j(t)}
\Gamma_{ij, ji} (t)dt
+ d \Upsilon_i(t).
\]
Here $d \Upsilon_i(t)$ denotes
the finite-variation part of 
$(U^{\dagger}(t) d M(t) U(t))_{ii}$, $t \geq 0$, 
$1 \leq i \leq N$.
\end{prop}

We will show two basic examples of $M(t) \in \sH_N, t \geq 0$
and applications of Proposition \ref{thm:Bru} \cite{KT04}.
Let $\nu \in \N_0 := \N \cup \{0\}$ and
$(B_{ij}(t))_{t \geq 0}$, 
$(\widetilde{B}_{ij}(t))_{t \geq 0}$, 
$1 \leq i \leq N+\nu$, $1 \leq j \leq N$
be independent one-dimensional standard Brownian motions.
For $1 \leq i \leq j \leq N$, put
\[
S_{ij}(t)= \begin{cases}
B_{ij}(t)/\sqrt{2}, \quad & (i<j),
\cr
B_{ii}(t), \quad & (i=j),
\end{cases}
\quad 
A_{ij}(t)= \begin{cases}
\widetilde{B}_{ij}(t)/\sqrt{2}, \quad & (i<j),
\cr
0, \quad & (i=j),
\end{cases}
\]
and let
$S_{ij}(t)=S_{ji}(t)$ and $A_{ij}(t)=-A_{ji}(t)$, $t \geq 0$
for $1 \leq j < i \leq N$.

\vskip 0.3cm
\noindent{\bf Example 2.1}
Put
$M_{ij}(t)=S_{ij}(t)+\sqrt{-1} A_{ij}(t), t \geq 0$, 
$1 \leq i, j \leq N$.
By definition,
$\langle d M_{ij}, d M_{k \ell} \rangle_t=\delta_{i \ell} \delta_{j k} dt$, $t \geq 0$, 
$1 \leq i, j, k, \ell \leq N$. Hence, by unitarity of $U(t), t \geq 0$, 
we see that $\Gamma_{ij, k \ell}(t) = \delta_{i \ell} \delta_{j k}$, 
which gives 
$\langle d \cM_i, d \cM_j \rangle_{t}=\Gamma_{ii, jj}(t) dt =\delta_{ij} dt$
and $\Gamma_{ij, ji}(t) \equiv 1$, $t \geq 0$, $1 \leq i, j \leq N$.
Then, Proposition \ref{thm:Bru} proves that
the eigenvalue process $(\bLambda(t))_{t \geq 0}$, 
satisfies the following system of SDEs with $\beta=2$,
\begin{equation}
d \Lambda_i(t)= d B_i(t) + \frac{\beta}{2} \sum_{1 \leq j \leq N, j \not=i}
\frac{dt}{\Lambda_i(t)-\Lambda_j(t)}, 
\label{eqn:Dyson1}
\end{equation}
$t \geq 0, 1 \leq i \leq N$.
Here, $(B_i(t))_{t \geq 0}, 1 \leq i \leq N$ are
independent one-dimensional standard Brownian motions,
which are different from $(B_{ij}(t))_{t\geq 0}$ 
and $(\widetilde{B}_{ij}(t))_{t \geq 0}$ used
to define $(S_{ij}(t))_{t \geq 0}$ and $(A_{ij}(t))_{t \geq0}$, 
$1 \leq i, j \leq N$.

\vskip 0.3cm
\noindent{\bf Example 2.2}
Consider an $(N+\nu) \times N$ rectangular-matrix-valued process 
given by 
$K(t)=(B_{ij}(t)
+\sqrt{-1} \widetilde{B}_{ij}(t))_{1 \leq i \leq N+\nu, 1 \leq j \leq N}$,
$t \geq 0$, 
and define an $\sH_N$-valued process by
$M(t)=K^{\dagger}(t) K(t), t \geq 0$.
The matrix $M(t)$ is positive semi-definite and hence 
the eigenvalues are non-negative;
$\Lambda_i(t) \in \R_+$, 
$t \geq 0, 1 \leq i \leq N$. 
We see that the finite-variation part of $dM_{ij}(t)$ is 
equal to $2(N+\nu) \delta_{ij} dt$, $t \geq 0$, 
and 
$\langle d M_{ij}, d M_{k \ell} \rangle_t=2(M_{i \ell}(t) 
\delta_{jk}+M_{k j}(t) \delta_{i \ell}) dt$,
$t \geq 0$, $1 \leq i, j, k, \ell \leq N$,
which implies that
$d \Upsilon_i(t)=2(N+\nu) dt$, 
$\Gamma_{ij, ji}(t)=2(\Lambda_i(t)+\Lambda_j(t))$, 
and $\langle d \cM_i, d \cM_j \rangle_t = \Gamma_{ii, jj}(t) dt 
= 4 \Lambda_i(t) \delta_{ij} dt$,
$t \geq 0$, $1 \leq i, j \leq N$.
Then, we have the SDEs for eigenvalue processes,
\begin{align}
d \Lambda_i(t)
&= 2 \sqrt{\Lambda_i(t)} d \widetilde{B}_i(t) 
\label{eqn:Wishart1}
\\
& \quad
+ \beta \Bigg[ (\nu+1) + 2 \Lambda_i(t)
\sum_{1 \leq j \leq N,  j \not=i}
\frac{1}{\Lambda_i(t)-\Lambda_j(t)} \Bigg] dt, 
\nonumber
\end{align}
$t \geq 0, 1 \leq i \leq N$ with $\beta=2$,
where $(\widetilde{B}_i(t))_{t \geq 0}, 1 \leq i \leq N$ are
independent one-dimensional standard Brownian motions,
which are different from $(B_{ij}(t))_{t \geq 0}$ 
and $(\widetilde{B}_{ij}(t))_{t \geq 0}$, 
$1 \leq i, j \leq N$, used above
to define the rectangular-matrix-valued process
$(K(t))_{t \geq 0}$.
The positive roots of eigenvalues of $M(t)$ give the
{\it singular values} of the rectangular matrix $K(t)$,
which are denoted by 
$\cS_i(t) =\sqrt{\Lambda_i(t)}, t \geq 0, 1 \leq i \leq N$.
The system of SDEs for them is readily obtained
from (\ref{eqn:Wishart1}) as
\begin{align}
d \cS_i(t) &= d \widetilde{B}_i(t)
+ \frac{\beta(\nu+1)-1}{2 \cS_i(t)} dt
\label{eqn:Wishart2}
\\
& \quad + \frac{\beta}{2} 
\sum_{1 \leq j \leq N, j \not=i}
\left( \frac{1}{\cS_i(t)-\cS_j(t)} + \frac{1}{\cS_i(t)+\cS_j(t)} \right) dt,
\nonumber
\end{align}
$t \geq 0, 1 \leq i \leq N$
with $\beta=2$ and $\nu \in \N_0$. 

Other examples of $\sH_N$-valued processes $(M(t))_{t \geq 0}$ are
shown in \cite{KT04}, in which the eigenvalue processes
following the SDEs (\ref{eqn:Dyson1}),
(\ref{eqn:Wishart1}), and (\ref{eqn:Wishart2}) 
with $\beta=1$ and 4 are
also shown.

\subsection{2D-Coulomb gases confined in 1D}
\label{sec:log_gases2}

In the next section, we will consider the SLE.
Schramm used a parameter $\kappa >0$ in order to parameterize
time changes of a Brownian motion \cite{Sch00}.
Accordingly, we relate the parameter $\beta$ to $\kappa$
by setting $\beta=8/\kappa$, 
and perform a time change $t \to \kappa t$.
Since 
$(B(\kappa t))_{t \geq 0} \law= (\sqrt{\kappa} B(t))_{t \geq 0}$,
if we put
$Y^{\R}_i(t) := \Lambda_i(\kappa t)$,
$Y^{\R_+}_i(t) := \cS_i(\kappa t)$, $t \geq 0$, $1 \leq i \leq N$,
the system of SDEs (\ref{eqn:Dyson1}) gives (\ref{eqn:Dyson0})
and that of (\ref{eqn:Wishart2})
gives
\begin{align}
& d Y^{\R_+}_i(t) = \sqrt{\kappa} d \widetilde{B}_i(t)
+ \frac{8(\nu+1)-\kappa}{2Y^{\R_+}_i(t)} dt
\label{eqn:Bru_WishartB}
\\
& \qquad + 4 \sum_{1 \leq j \leq N, j \not=i}
\Bigg( \frac{1}{Y^{\R_+}_i(t)-Y^{\R_+}_j(t)}
+ \frac{1}{Y^{\R_+}_i(t)+Y^{\R_+}_j(t)} \Bigg) dt, 
\nonumber
\end{align}
$t \geq 0, 1 \leq i \leq N$, where $\nu \geq 0$.
In the present paper, we call
$(\Y^{\R}(t))_{t \geq 0}$
the $(8/\kappa)$-{\it Dyson model} and
$(\Y^{\R_+}(t))_{t \geq 0}$
the $(8/\kappa, \nu)$-{\it Bru--Wishart process},
respectively.
The above systems of SDEs for
$(\Y^{S}(t))_{t \geq 0}$
can be written as
\[
dY^S_i(t)=\sqrt{\kappa} dB_i(t)
+ \left. 
\frac{\partial \phi^S(\x)}{\partial x_i}
\right|_{\x=\Y^S(t)}dt, \quad
t \geq 0, \quad 1 \leq i \leq N, 
\]
$S=\R$ or $\R_+$, 
when we introduce the following logarithmic potentials,
\begin{equation}
\phi^{S}(\x) :=
\begin{cases}
\displaystyle{4 \sum_{1 \leq i < j \leq N} 
\log(x_j-x_i)},
\qquad \mbox{for $S=\R$}, \\
\displaystyle{4 \sum_{1 \leq i < j \leq N}
\Big[ \log(x_j-x_i)+\log(x_j+x_i) \Big]} \\[2mm] 
\quad + \displaystyle{\frac{8(\nu+1)-\kappa}{2}} 
\sum_{1 \leq i \leq N} \log x_i, 
\quad \mbox{for $S=\R_+$}.
&
\end{cases}
\label{eqn:potential}
\end{equation}

In this sense,
the $(8/\kappa)$-Dyson model
and the $(8/\kappa, \nu)$-Bru--Wishart process 
are regarded as {\it stochastic log-gases} 
in one dimension \cite{For10}.
Since the logarithmic potential describes
the two-dimensional Coulomb law in electrostatics,
the present processes are also considered
as stochastic models of 
{\it 2D-Coulomb gases confined in 1D}.

\SSC{Multiple Schramm--Loewner Evolution} \label{sec:multiple_SLE}
\subsection{Loewner equations for a single-slit and a multi-slit}
\label{sec:chordal_LE}

Let $D$ be a simply connected domain $D \subsetneq \C$
with boundary $\partial D$.
We consider a slit in $D$, which is defined as a simple curve 
$\eta=\{\eta(t) : t \in (0, \infty)\} \subset D$; 
$\eta(s) \not= \eta(t)$ for $s \not=t$
and suppose that $\lim_{t \to 0} \eta(t)=: \eta(0) \in \partial D$.
Let $\eta(0,t]:=\{\eta(s) : s \in (0, t]\}$ and 
$D^{\eta}_t := D \setminus \eta(0, t], t \in (0, \infty)$
with $D^{\eta}_0 := D$. 
The Loewner theory describes the slit $\eta$
by encoding it into a time-dependent
analytic function $(g_{D^{\eta}_t})_{t \geq 0}$
such that 
\[
g_{D^{\eta}_t} : \mbox{conformal map} \, \, 
D^{\eta}_t \to D,
\quad t \in [0, \infty).
\]
Let us apply the Loewner theory to the case of $D=\H$, 
in which $\eta(0) \in \R$
and $\eta \subset \H$.
Let $\H^{\eta}_t := \H \setminus \eta(0, t]$, $t > 0$
and $\H^{\eta}_0 := \H$.
Then, for each time  $t \geq 0$,
$\H^{\eta}_t$ is a simply connected
domain in $\C$ and there exists a unique conformal map
$\H^{\eta}_t \to \H$ 
satisfying the condition
$g_{\H^{\eta}_t}(z)=z+ \hcap(\eta(0, t])/z + \rO(|z|^{-2})$
as $z \to \infty$, $t > 0$,
in which the coefficient of $z$ is unity and 
no constant term appears.
This is called the
{\it hydrodynamic normalization}
and $\hcap(\eta(0, t])$ gives
the {\it half-plane capacity} of 
$\eta(0, t]$.
The following can be proved 
(see, for instance, \cite[Proposition 4.4]{Law05}).

\begin{thm}
\label{thm:LE}
Let $\eta$ be a slit in $\H$ such that
$\hcap(\eta(0, t])=2t$, $t >0$.
Then, the solution $(g_t)_{t \geq 0}$ of the differential equation
(chordal Loewner equation) 
\begin{equation}
\frac{d g_t(z)}{dt} = \frac{2}{g_t(z)-V(t)}, \quad
t \geq 0, \quad g_0(z)=z,
\label{eqn:LE1}
\end{equation}
where
\[
V(t)=g_{\H^{\eta}_t}(\eta(t))
:= \lim_{z \to \eta(t), z \in \H^{\eta}_t}
g_{\H^{\eta}_t}(z), \quad t \geq 0,
\]
coincides with $(g_{\H^{\eta}_t})_{t \geq 0}$.
\end{thm}

Theorem \ref{thm:LE} can be extended 
to the situation such that $\eta$ is given by
a multi-slit $\bigcup_{i=1}^N \eta_i \subset \H$ and
$\H^{\eta}_t := \H \setminus \bigcup_{i=1}^N \eta_i(0, t]$,
$t > 0$ with $\H^{\eta}_0 :=\H$. 

\begin{thm}
\label{thm:mLE}
For $N \in \N$, 
let $\bigcup_{i=1}^N \eta_i$ be a multi-slit in $\H$ such that
$\hcap(\bigcup_{i=1}^N \eta(0, t])=2t, t >0$.
Then, there exists a set of weight functions 
$(\lambda_i(t))_{t \geq 0}$, $1 \leq i \leq N$
satisfying 
$\lambda_i(t) \geq 0, 1 \leq i \leq N$,
$\sum_{i=1}^N \lambda_i(t)=1, t \geq 0$
such that the 
solution $(g_t)_{t \geq 0}$ of the differential equation
(multiple chordal Loewner equation) 
\begin{equation}
\frac{d g_t(z)}{dt} = 
\sum_{i=1}^N \frac{2 \lambda_i(t)}{g_t(z)-V_i(t)}, \quad
t \geq 0, \quad 
g_0(z)=z,
\label{eqn:mLE1}
\end{equation}
where 
\[
V_i(t)=g_{\H^{\eta}_t}(\eta_i(t))
:= \lim_{z \to \eta_i(t), z \in \H^{\eta}_t}
g_{\H^{\eta}_t}(z), 
\quad t \geq 0, \quad 1 \leq i \leq N,
\]
coincides with $(g_{\H^{\eta}_t})_{t \geq 0}$. 
\end{thm}
\noindent{\it Proof} \, \,
We parameterize each curve $\eta_{i}$ separately 
by $\eta_{i}: (0,\infty)\to\H; s_{i}\mapsto \eta_{i}(s_{i})$ 
so that $\hcap \Big(\bigcup_{j=1}^{N}\eta_{j}(0,s_{j}] \Big)$ is 
differentiable with respect to $s_{i}$, $1 \leq i \leq N$. 
For each $\Us=(s_{1},\dots, s_{N}) \in [0, \infty)^{N}$, 
we set 
$\H_{\Us}^{\eta}=\H \backslash \bigcup_{i=1}^{N}\eta_{i}(0,s_{i}]$. 
A similar argument as in \cite[Proposition 4.4]{Law05} shows 
that the family of conformal maps 
$g_{\Us}=g_{\H^{\eta}_{\Us}}$, $\Us\in [0, \infty)^{N}$
satisfies partial differential equations
\[
\frac{\partial g_{\Us}(z)}{\partial s_{i}}
=\frac{1}{g_{\Us}(z)-V_{i}(\Us)}
\frac{\partial}{\partial s_{i}}
\hcap 
\Bigg(\bigcup_{j=1}^{N}\eta_{j}(0,s_{j}] \Bigg),
\quad 1 \leq i \leq N, 
\]
where
\[
V_{i}(\Us)=g_{\H^{\eta}_{\Us}}(\eta_i(s_{i}))
:= \lim_{z \to \eta_i(s_{i}), z \in \H^{\eta}_{\Us}}
g_{\H^{\eta}_{\Us}}(z).
\]
We have these parameters dependent 
on a single parameter so that 
$s_{i}=s_{i}(t)$, $t\geq 0$, $1 \leq i \leq N$ 
are increasing and differentiable in $t$, 
and write $\Us=\Us(t)$, $t\geq 0$. 
Then, we may understand $\H^{\eta}_{t}=\H^{\eta}_{\Us(t)}$, 
$g_{t}=g_{\Us(t)}$, and $V_{i}(t)=V_{i}(\Us(t))$, 
$1 \leq i \leq N$, $t\geq 0$. 
Furthermore, we impose a condition that 
$\hcap \Big(\bigcup_{j=1}^{N} \eta_{j}(0,s_{j}(t)] \Big)=2t$. 
Then, the family of conformal maps $(g_{t})_{t \geq 0}$ 
satisfies the desired differential equation (\ref{eqn:mLE1}), 
where
\[
\lambda_{i}(t)=\frac{1}{2} 
\frac{\partial}
{\partial s_{i}} 
\hcap \Bigg(\bigcup_{j=1}^{N}\eta_{j}(0,s_{j}] \Bigg)
\Bigg|_{\Us=\Us (t)}
\frac{ds_{i}(t)}{dt}, \ \ 1 \leq i \leq N, \ \ t\geq 0
\]
are subject to the constraint $\sum_{i=1}^{N}\lambda_{i}(t)=1$, 
$t\geq 0$.
\qed

The multiple chordal Loewner equation (\ref{eqn:mLE1}) 
for $D=\H$ can be mapped to other simply 
connected proper subdomains of $\C$ by conformal maps.
Here, we consider a conformal map
$\widehat{\varphi}(z)=z^2 : \O \to \H$.
We set 
$\widehat{g}_t(z) =\sqrt{g_t(z^2)+c(t)}$, $t \geq 0$ 
with a function of time $c(t), t \geq 0$.
Then, we can see that (\ref{eqn:mLE1}) is transformed to
\begin{equation}
\frac{d \widehat{g}_t(z)}{dt}
=\sum_{i=1}^N \left( 
\frac{2 \widehat{\lambda}_i(t)}{\widehat{g}_t(z)-\widehat{V}_i(t)}
+\frac{2 \widehat{\lambda}_i(t)}{\widehat{g}_t(z)+\widehat{V}_i(t)} 
\right)
+\frac{2 \widehat{\lambda}_0(t)}{\widehat{g}_t(z)}, \, t \geq 0,
\label{eqn:LE_O1}
\end{equation}
$\widehat{g}_0(z) =z \in \O$, 
where 
$\widehat{V}_i(t)=\sqrt{V_i(t)+c(t)}$, 
$t \geq 0, 1 \leq i \leq N$ 
and 
$2 \sum_{i=1}^N \widehat{\lambda}_i(t)
+\widehat{\lambda}_0(t)=(1/4) dc(t)/dt, t \geq 0$.
Here, we can assume that $\widehat{V}_i(t) \in \R_+$
by a proper choice of the function $c(t)$, $t\geq 0$.
The equation (\ref{eqn:LE_O1})
can be regarded as a
multi-slit version of the
{\it quadrant Loewner equation} studied in \cite{Tak14}.
The solution of (\ref{eqn:LE_O1}) gives 
a conformal map 
$\widehat{g}_t=g_{\O^{\eta}_t} : 
\O^{\eta}_t \to \O$,
where $\O^{\eta}_t := \O \setminus \bigcup_{i=1}^N \eta_i(0, t]$, $t > 0$,
$\O^{\eta}_0 :=\O$, and 
$g_{\O^{\eta}_t}(\eta_i(t))=\widehat{V}_i(t) \in \R_+$,
$t \geq 0$, $1 \leq i \leq N$. 

\subsection{SLE}
\label{sec:SLE}

So far we have considered the problem in which 
given a single slit $\eta(0, t], t > 0$
or a multi-slit $\bigcup_{i=1}^N \eta(0, t], t >0$ in $\H$,
the time-evolution of the conformal map from $\H^{\eta}_t$ to $\H$,
$t \geq 0$ is asked.
The answers are given by Theorem \ref{thm:LE} 
and Theorem \ref{thm:mLE}.
For $\H$ with a single slit, Schramm considered
the inverse problem in a probabilistic setting \cite{Sch00}.
He first asked a suitable family of driving stochastic processes
$(X(t))_{t \geq 0}$ on $\R$.
Then, he asked the probability law of 
the random slit $\eta$ in $\H$
that is determined by the solution $g_t=g_{\H^{\eta}_t}, t \geq 0$
of the Loewner equation (\ref{eqn:LE1}) via 
$X(t)=g_{\H^{\eta}_t}(\eta(t)), t \geq 0$.
Schramm argued that the conformal invariance 
and the domain Markov property 
of the law of the curve imply that the
driving process $(X(t))_{t \geq 0}$ should be 
$(B(\kappa t))_{t \geq 0} \law= (\sqrt{\kappa} B(t))_{t \geq 0}$
with a parameter $\kappa >0$.
The solution of the chordal Loewner equation (\ref{eqn:LE1}) 
driven by $X(t)=\sqrt{\kappa} B(t), t \geq 0$
is called the {\it chordal Schramm--Loewner evolution} 
with parameter $\kappa >0$ and is written as
chordal SLE($\kappa$) for short. 

The following was proved by Lawler, Schramm, 
and Werner \cite{LSW04}
for $\kappa=8$ and by Rohde and Schramm 
\cite{RS05} for $\kappa \not=8$.
\begin{prop} 
\label{thm:SLE_curve}
A chordal SLE($\kappa$)
$(g_{\H^{\eta}_t})_{t \geq 0}$ determines 
a continuous curve 
$\eta=\{\eta(t): t \in [0, \infty)\} \subset \overline{\H}$
a.s.
\end{prop}

In this inverse problem, the domain $\H^{\eta}_{t}$ is defined 
as the unbounded component of 
$\H\backslash \eta(0,t]$ so that the solution 
$g_{t}$ is a conformal map from 
$\H^{\eta}_{t}$ to $\H$ at each $t>0$, which verifies 
writing the solution as $g_{\H^{\eta}_{t}}=g_{t}$, $t\geq 0$.

The continuous curve $\eta$ determined 
by an SLE($\kappa$) 
is called an {\it SLE($\kappa$) curve}.
The probability law of an SLE($\kappa$) curve qualitatively 
depends on $\kappa$. 
When $\kappa \in (0, 4]$, the SLE($\kappa$) curve is
a simple curve in $\H$.
It becomes self-intersecting and can touch the real axis $\R$
when $\kappa > 4$, 
and becomes a space-filling curve when $\kappa \geq 8$
(see, for instance, \cite{Law05,Kat15_Springer}).

\subsection{Multiple SLE}
\label{sec:mSLE}

For simplicity, we assume that
$\lambda_i(t) \equiv 1/N, t \geq 0, 1 \leq i \leq N$
in (\ref{eqn:mLE1}) in Theorem \ref{thm:mLE}.
Then, by a simple time change $t/N \to t$
associated with a change of notation, $g_{Nt} \to g_{\H^{\eta}_t}$,
the multiple chordal Loewner equation is written as (\ref{eqn:mSLE0}). 
Then, we ask what is a suitable family of driving
stochastic processes of $N$ particles
$(\X^{\R}(t))_{t \geq 0}$ on $\R$.

Bauer, Bernard, and Kyt\"ol\"a \cite{BBK05}
and Graham \cite{Gra07} argued that 
$X^{\R}_i(t)$, $t \geq 0$, $1 \leq i \leq N$
are semi-martingales and the quadratic variations 
should be given by 
$\langle d X^{\R}_i, dX^{\R}_j \rangle_t= \kappa \delta_{ij} dt, t \geq 0$,
$1 \leq i, j \leq N$ with $\kappa > 0$.
Then, we can assume that the system of SDEs for
$(\X^{\R}(t))_{t \geq 0}$ is in the form (\ref{eqn:SDE0}). 

In the orthant system (\ref{eqn:LE_O1}),
we put $\widehat{\lambda}_i(t) \equiv r/(2N)$, 
$t \geq 0, r \in (0, 1], 1 \leq i \leq N$, $dc(t)/dt \equiv 4, t \geq 0$, 
and perform a time change
$r t/(2N) \to t$ associated with a change of
notation $\widehat{g}_{2N t/r} \to g_{\O^{\eta}_t}$.
Then, the multiple Loewner equation in $\O$ is written as
\begin{equation}
\frac{d g_{\O^{\eta}_t}(z)}{dt}
= \sum_{i=1}^N \Bigg( 
\frac{2}{g_{\O^{\eta}_t}(z)-X^{\R_+}_i(t)} 
+ \frac{2}{g_{\O^{\eta}_t}(z)+X^{\R_+}_i(t)} \Bigg) 
+ \frac{4 \delta}{g_{\O^{\eta}_t}(z)}, 
\label{eqn:mSLE2}
\end{equation}
$t \geq 0$ with $g_{\O^{\eta}_0}(z) =z \in \O$, 
where $\delta := N(1-r)/r \geq 0$.
We assume that the system of SDEs for
$\X^{\R_+}(t) \in (\R_+)^N, t \geq 0$ 
is in the form (\ref{eqn:SDE0}).

Analogously to the case of the SLE for a single slit,
in both cases of $D=\H$ and $\O$, we find a family of domains $(D^{\eta}_{t}\subseteq D)_{t\geq 0}$ such that
$g_{D^{\eta}_{t}}$ is a conformal map from 
$D^{\eta}_{t}$ to $D$ at each $t\ge 0$.
(See Remark 1.4 (1) in \cite{KK21}.)

\SSC{Gaussian Free Field with Zero Boundary Condition } \label{sec:GFF}
\subsection{Bochner--Minlos Theorem} 
\label{sec:GFF_Dirichlet}

Let $D \subsetneq \C$ be a simply connected domain.
Consider the real $L^2$ space  
with the inner product, 
$(f, g) := \int_D f(z) g (z) d \mu(z)$, $f, g \in L^2(D)$,
where $\mu(z)$ is the Lebesgue measure on $\C$;
$d \mu(z)=\sqrt{-1} dz d\overline{z}/2$.
Let $\Delta$ be the Dirichlet Laplacian acting on $L^2(D)$.
In the present subsection \ref{sec:GFF_Dirichlet}
we assume that $D$ is bounded.
Then $-\Delta$ has positive discrete eigenvalues so that
$-\Delta e_n = \lambda_n e_n$, 
$e_n \in L^2(D)$, $n \in \N$.
We assume that the eigenvalues are labeled 
in the non-decreasing order;
$0 < \lambda_1 \leq \lambda_2 \leq \cdots$.
The system of eigenfunctions $\{e_n\}_{n \in \N}$ forms
a CONS of $L^2(D)$. 
The asymptotic behavior of eigenvalues obeys
{\it Weyl's formula}; $\lim_{n \to \infty} \lambda_n/n =\rO(1)$.

For $f, g \in \cC_{\rm c}^{\infty}(D)$,
the {\it Dirichlet inner product} is defined by
\begin{equation}
(f, g)_{\nabla} := \frac{1}{2 \pi} \int_D (\nabla f)(z) \cdot (\nabla g)(z)
d \mu(z).
\label{eqn:Dirichlet_IP}
\end{equation}
The Hilbert space completion of $\cC_{\rm c}^{\infty}(D)$ with respect
to $(\cdot, \cdot)_{\nabla}$ will be denoted by $W(D)$.
We write 
$\|f\|_{\nabla}=\sqrt{(f, f)_{\nabla}}, f \in W(D)$.
If we set $u_n=\sqrt{2\pi/\lambda_n} \, e_n, n \in \N$,
then, by integration by parts, we have
$(u_n, u_n)_{\nabla}=
(u_n, (-\Delta) u_m)/(2\pi)=\delta_{nm}$,
$n, m \in \N$.
Therefore, $\{u_n\}_{n \in \N}$ forms a CONS of $W(D)$.

Let $\widehat{\cH}(D)$ be the space of formal infinite
series in $\{u_n\}_{n \in \N}$, 
which is obviously isomorphic to $\R^{\N}$
by setting 
$\widehat{\cH}(D) \ni \sum_{n \in \N} f_n u_n 
\mapsto (f_n)_{n \in \N} \in \R^{\N}$.
As a subspace of $\widehat{\cH}(D)$, $W(D)$ is isomorphic to
$\ell^2(\N) \subset \R^{\N}$.
For two formal series 
$f=\sum_{n \in \N} f_n u_n$,
$g=\sum_{n \in \N} g_n u_n \in \widehat{\cH}(D)$
such that
$\sum_{n \in \N} |f_n g_n| < \infty$,
we define their pairing as
$(f, g)_{\nabla}:=\sum_{n \in \N} f_n g_n$.
In the case when $f, g \in W(D)$, their pairing, of course, 
coincides with the Dirichlet inner product 
(\ref{eqn:Dirichlet_IP}).

Notice that, for any $a \in \R$, the operator
$(-\Delta)^a$ acts on $\widehat{\cH}(D)$ as
$(-\Delta)^a \sum_{n \in \N} f_n u_n
:= \sum_{n \in \N} \lambda_n^a f_n u_n$,
$(f_n)_{n \in \N} \in \R^{\N}$.
Using this fact, we define
$\cH_a(D) := (-\Delta)^a W(D)$, $a \in \R$,
each of which is a Hilbert space with the inner product
$\langle f, g \rangle_a
:=((-\Delta)^{-a}f, (-\Delta)^{-a}g)_{\nabla}$,
$f, g \in \cH_a(D)$.
We can prove that $\cH_a(D) \subset \cH_b(D)$ for $a<b$ 
using Weyl's formula for $\{\lambda_n\}_{n \in \N}$,
and that
the dual Hilbert space of $\cH_a(D)$ is given by $\cH_{-a}(D)$
(see \cite{Ara10}). 

\vskip 0.3cm
\noindent{\bf Remark 4.1} 
Since
\[ 
\langle f, g \rangle_{1/2}
=( (-\Delta)^{-1/2} f, (-\Delta)^{-1/2}g)_{\nabla}
=(f, g)/(2 \pi), 
\quad 
f, g \in \cH_{1/2}(D),
\]
$\cH_{1/2}(D) \simeq L^2(D)$.
This implies that  
the members of $\cH_a(D)$
with $a > 1/2$ cannot be functions, 
but are distributions.
\vskip 0.3cm

Define $\cE(D) := \bigcup_{a > 1/2} \cH_a(D)$.
Then, its dual Hilbert space is identified with
$\cE(D)^{\ast} :=\bigcap_{a < -1/2} \cH_a(D)$ and 
$\cE(D)^{\ast} \subset W(D) \subset \cE(D)$
is established.
Here $(\cE(D)^{\ast}, W(D), \cE(D))$ is called a {\it Gel'fand triple}.
We set 
$\Sigma_{\cE(D)} =\sigma(\{(\cdot, f)_{\nabla} : f \in \cE(D)^{\ast} \})$.
On such a setting, the following is proved.
This theorem is called 
the {\it Bochner--Minlos theorem} \cite{Hid80,She07,Ara10}. 

\begin{thm}[Bochner--Minlos theorem]
\label{thm:BM}
Let $\psi$ be a continuous function of positive type
on $W(D)$ such that $\psi(0)=1$.
Then there exists a unique probability measure 
$\bP$ on $(\cE(D), \Sigma_{\cE(D)})$
such that
$\psi(f)=\int_{\cE(D)} e^{\sqrt{-1} (h, f)_{\nabla}} \bP(d h)$
for $f \in \cE(D)^{\ast}$.
\end{thm}

Under certain conditions on $\psi$,
the domain of the random functional $h$ in the above formula
can be extended from $\cE(D)^{\ast}$ to $W(D)$.
It is easy to verify that the functional
$\Psi(f):= e^{-\|f\|_{\nabla}^2/2}$ satisfies the conditions.
Then, the following is established
with the probability measure $\bP$ on $(\cE(D), \Sigma_{\cE(D)})$,
\begin{equation}
\int_{\cE(D)} e^{\sqrt{-1} (h, f)_{\nabla}} \bP(d h)
=e^{-\|f\|_{\nabla}^2/2}
\quad \mbox{for $f \in W(D)$}.
\label{eqn:BMeq}
\end{equation}

\begin{df}[zero-boundary GFF]
\label{thm:GFF_Dirichlet}
A Gaussian free field (GFF) with zero
boundary condition (zero-boundary GFF) is defined as a pair 
$((\Omega^{\GFF}, \cF^{\GFF}, \P^{\GFF}), H)$ 
of a probability space 
$(\Omega^{\GFF},\cF^{\GFF}, \P^{\GFF})$ 
and an isometry
$H : W(D) \to L^{2}(\Omega^{\GFF}, \cF^{\GFF}, \P^{\GFF})$ 
such that each $H(f)$, $f\in W(D)$ is a centered Gaussian random variable.
\end{df}

For each $f \in W(D)$, we write 
$(H,f)_{\nabla} \in L^{2}(\cE(D),\Sigma_{\cE(D)}, \bP)$ 
for the random variable defined by 
$h \mapsto (h,f)_{\nabla}$, $h \in \cE(D)$.
Then (\ref{eqn:BMeq}) ensures that the pair of 
$((\cE(D),\Sigma_{\cE(D)}, \bP),H)$
gives a GFF with zero boundary condition. 
We often just call $H$ 
a zero-boundary GFF without referring 
to the probability space 
$(\Omega^{\GFF}, \cF^{\GFF}, \P^{\GFF})
=(\cE(D),\Sigma_{\cE(D)}, \bP)$.

\subsection{Conformal invariance of a zero-boundary GFF} 
\label{sec:conformal_invariance}

Assume that $D, \widetilde{D} \subsetneq \C$ are simply connected domains
and let $\varphi: \widetilde{D} \to D$ be a conformal map.

\begin{lem}
\label{thm:conformal_inv}
The Dirichlet inner product (\ref{eqn:Dirichlet_IP}) is conformally invariant.
That is, for $f, g \in \cC_{\rm c}^{\infty}(D)$, 
\[
\int_{D} (\nabla f)(z) \cdot (\nabla g)(z) d \mu(z)
=\int_{\widetilde{D}} (\nabla (f \circ \varphi))(z) \cdot 
(\nabla (g \circ \varphi))(z) d \mu(z). 
\]
\end{lem}

From the above lemma, we see that
$\varphi^{\ast} : W(D) \ni f \mapsto f \circ \varphi \in W(\widetilde{D})$
is an isomorphism. This allows one to consider
a GFF on an unbounded domain.
Namely, if $\widetilde{D}$ is bounded on which 
a zero-boundary GFF is defined, 
but $D$ is unbounded, we can define a family
$\{(\varphi_{\ast} H, f)_{\nabla} : f \in W(D)\}$ by
$(\varphi_{\ast} H, f)_{\nabla} := (H, \varphi^{\ast} f)_{\nabla}, f \in W(D)$
so that we have the covariance structure,
\[ 
\E^{\GFF} \Big[
(\varphi_{\ast} H, f)_{\nabla} (\varphi_{\ast} H, g)_{\nabla} \Big]
=(\varphi^{\ast} f, \varphi^{\ast} g)_{\nabla}
=(f, g)_{\nabla}, \quad f, g \in W(D),
\]
where $\E^{\GFF}$ is the expectation value 
with respect to $\P^{\GFF}$.
Relying on the formal computation, 
\begin{align*}
(\varphi_{\ast} H, f)_{\nabla}
=(H, \varphi^{\ast} f)_{\nabla}
&=\frac{1}{2 \pi} \int_{\widetilde{D}} (\nabla H)(z) \cdot (\nabla f \circ \varphi)(z) d \mu(z)
\nonumber\\
&=\frac{1}{2 \pi} \int_D (\nabla H \circ \varphi^{-1})(z) \cdot (\nabla f)(z) d \mu(z), 
\end{align*}
we understand the equality $\varphi_{\ast} H = H \circ \varphi^{-1}$.
By the fact shown above that
the covariance structure does not change under 
a conformal map $\varphi$, we say 
{\it a zero-boundary GFF is conformally invariant}. 

\subsection{Green's function of a zero-boundary GFF} 
\label{sec:Green}

Assume that $D \subsetneq \C$ is a simply connected domain.
In the previous subsections, we have constructed a family
$\{(H, f)_{\nabla} : f \in W(D)\}$ of random variables
whose covariance structure is given by
\[ 
\E^{\GFF} \Big[
(H, f)_{\nabla} (H, g)_{\nabla} \Big]=(f, g)_{\nabla}, 
\quad 
f, g \in W(D).
\]
By a formal integration by parts, we see that
\begin{align*}
(H, f)_{\nabla} &= \frac{1}{2\pi} \int_{D} (\nabla H)(z) 
\cdot (\nabla f)(z) d \mu(z)
=\frac{1}{2\pi} \int_{D} H(z) (-\Delta f)(z) d \mu(z)
\nonumber\\
&=\frac{1}{2 \pi} (H, (-\Delta) f).
\end{align*}
Motivated by this observation, we define
\begin{equation}
(H, f):=2 \pi (H, (-\Delta)^{-1} f)_{\nabla} 
\quad \mbox{for $f \in \sD((-\Delta)^{-1})$}, 
\label{eqn:inner_product}
\end{equation}
where
$\sD((-\Delta)^{-1})$ denotes the domain of
$(-\Delta)^{-1}$ in $W(D)$.
The action of $(-\Delta)^{-1}$ is expressed as
an integral operator as
\[ 
((-\Delta)^{-1} f)(z)
=\frac{1}{2 \pi} \int_{D} G_D(z,w) f(w) d \mu(w)
\quad \mbox{a.e. $z \in D$}, 
\]
$f \in \sD((-\Delta)^{-1})$,
where the integral kernel
$G_D$ is known as {\it the Green's function}
of $D$ under the Dirichlet boundary condition:
$G_D(z, w)=0, w \in D$ if $z \in \partial D$.
Hence the covariance of $(H, f)$ and $(H, g)$ with
$f, g \in \sD((-\Delta)^{-1})$ is written as
\begin{equation}
\E^{\GFF}[ (H, f) (H, g)] 
=\int_{D \times D} f(z) G_D(z, w) g(w) d \mu(z) d \mu(w).
\label{eqn:Green}
\end{equation}
When we symbolically write
\[
(H, f)=\int_{D} H(z) f(z) d \mu(z),
\quad f \in \sD((-\Delta)^{-1}),
\]
the covariance structure can be expressed as
\[
\E^{\GFF}[H(z) H(w)]= G_D(z, w), \quad
z, w \in D, \quad z \not=w.
\]
The conformal invariance of a
zero-boundary GFF implies that
for a conformal map
$\varphi: \widetilde{D} \to D$, we have the equality,
\[
G_{\widetilde{D}}(z, w)=G_D(\varphi(z), \varphi(w)),
\quad z, w \in \widetilde{D}, 
\quad z \not= w.
\]

\vskip 0.3cm
\noindent{\bf Example 4.1} \,
When $D=\H$,
\[
G_{\H}(z, w)=
\log \left| \frac{z-\overline{w}}{z-w}
\right|,
\quad z, w \in \H, \quad z \not= w.
\]

\vskip 0.3cm
\noindent{\bf Example 4.2} \,
When $D=\O$,
\[
G_{\O}(z, w)=
\log \left|
\frac{(z-\overline{w})(z+\overline{w})}
{(z-w)(z+w)} 
\right|,
\quad z, w \in \O, \quad z \not= w.
\] 
\vskip 0.3cm

From the formula (\ref{eqn:Green}), we see that
$\cC_{\rm c}^{\infty}(D) \subset \sD((-\Delta)^{-1})$.
In the following, we will consider the
family of random variables $\{(H, f) : f \in \cC_{\rm c}^{\infty}(D) \}$
to characterize a GFF $H$.

\SSC{Gaussian Free Fields Coupled with Multiple SLEs} 
\label{sec:GFF_loggas}

In this section, we take a probability space 
$(\Omega, \cF, \P)$ on which a zero-boundary GFF $H$ 
and a multiple SLE $(g_{D^{\eta}_{t}})_{t \geq 0}$ are defined 
in such a way that they are independent, and 
a filtration $(\cF_{t})_{t\geq 0}$ 
to which the multiple SLE  is 
adapted 
(see Section \ref{sec:introduction} for a precise setting). 
We fix an open set $A\subset D$, where $D=\H$ or $\O$. 
Then the $(\cF_{t})_{t\geq 0}$-stopping time $\tau_{A}$ 
is defined in exactly the same expression
as (\ref{eqn:tau_A_multiple}). 
We also take an $(\cF_{t})_{t\geq 0}$-stopping time $\tau$ 
such that $\tau\le \tau_{A}$ a.s. as in Theorem \ref{thm:main0}. 

\subsection{Zero-boundary GFF transformed by a multiple SLE}
\label{sec:GFF_SLE}

Here, we write the zero-boundary GFF defined 
on $D=\H$ or $\O$ as $H_D$.
Consider the transformation of $H_D$ by the multiple SLE,
$H_{D^{\eta}_t} := H_D \circ g_{D^{\eta}_t}, 0\leq t\leq \tau$ 
on $D^{\eta}_{\tau}$.
By the conformal invariance, the Green's function of 
$H_{D^{\eta}_t}, 0\leq t\leq \tau$
is given by
$G_{D^{\eta}_t}(z, w) = G_D(g_{D^{\eta}_t}(z), g_{D^{\eta}_t}(w))$,
$z, w \in D^{\eta}_{\tau}$, $z \not=w$, $0\leq t\leq \tau$.
The following is obtained.

\begin{lem}
\label{thm:dG}
For $D=\H$ and $\O$, the increments of $G_{D^{\eta}_t}(z, w)$,
$z, w \in A$ in time $0\leq t\leq \tau$ are given as
\begin{align*}
d G_{\H^{\eta}_t}(z, w)
&=- \sum_{i=1}^N \Im \frac{2}{g_{\H^{\eta}_t}(z)-X^{\R}_i(t)}
\Im \frac{2}{g_{\H^{\eta}_t}(w)-X^{\R}_i(t)} dt,
\nonumber\\
d G_{\O^{\eta}_t}(z, w)
&=- \sum_{i=1}^N \Im 
\Bigg( \frac{2}{g_{\O^{\eta}_t}(z)-X^{\R_+}_i(t)}
- \frac{2}{g_{\O^{\eta}_t}(z)+X^{\R_+}_i(t)} \Bigg)
\nonumber\\
& \qquad \quad \times
\Im \Bigg( \frac{2}{g_{\O^{\eta}_t}(w)-X^{\R_+}_i(t)}
- \frac{2}{g_{\O^{\eta}_t}(w)+X^{\R_+}_i(t)} \Bigg) dt.
\end{align*}
\end{lem}
\noindent{\it Proof} \,
Using the explicit expressions of the Green's functions
given in Examples 4.1 and 4.2 
and the multiple Loewner equations (\ref{eqn:mSLE0}) 
and (\ref{eqn:mSLE2}), 
the increments of 
$(G_{D^{\eta}_t})_{0\leq t\leq \tau}$ 
with $D=\H$ and $\O$ are calculated.
The above expressions are obtained using the equality 
$\Re \zeta \overline{\omega}-\Re \zeta \omega
=2 \Im \zeta \Im \omega$ for $\zeta, \omega \in \C$.
\qed

\subsection{$\C$-valued logarithmic potentials and martingales}
\label{sec:martingales}

We have remarked in Section \ref{sec:log_gases2} that
the Dyson model and the Bru--Wishart process
studied in random matrix theory can be regarded
as stochastic log-gasses defined on a line $S=\R$
and a half-line $S=\R_+$, respectively.
There, the logarithmic potentials are given by (\ref{eqn:potential}). 
Here, we consider a {\it complex-valued logarithmic potentials}
between a point $z$ in the two-dimensional domain
$D \subsetneq \C$ and $N$ points 
$\x=(x_1, \dots, x_N)$ on the boundary $S$.
For $(D, S)=(\H, \R)$ and 
$(\O, \R_+)$, 
put
\begin{align*}
\Phi_{\H}(z, \x) &= \sum_{i=1}^N \log(z-x_i),
\nonumber\\
\Phi_{\O}(z, \x) &=
\Phi_{\O}(z, \x; q)
=\sum_{i=1}^N \Big\{ \log(z-x_i)+\log(z+x_i) \Big\}
+q \log z, 
\end{align*}
where $z \in D, \x \in S^N$, and $q \in \R$.

Now, we consider a time evolution of the $\C$-valued potential
$\Phi_{D}$ by letting $\x$ be the driving process
$(\X^S(t))_{t \geq 0}$ of the multiple SLE
$(g_{D^{\eta}_t})_{t \geq 0}$ and by transforming the
function $\Phi_{D}(\cdot, \X^{S}(t))$ by 
$(g_{D^{\eta}_t})_{t \geq 0}$.
We obtain the following.
\begin{lem}
\label{thm:c-potential_time}
For $D=\H$ and $\O$, 
the increments of the $\C$-valued potentials are
given as follows.
For $z \in A$, $\X^S(t) \in \W_N(S)$, 
$0\leq t\leq \tau$, 
\begin{align}
&
d \Phi_{\H}(g_{\H^{\eta}_t}(z), \X^{\R}(t))
=- \sum_{i=1}^N 
\frac{\sqrt{\kappa} dB_i(t)}{g_{\H^{\eta}_t}(z)-X^{\R}_i(t)} 
-\left(1 - \frac{\kappa}{4} \right) d \log g_{\H^{\eta}_t}'(z)
\label{eqn:dPhiH} 
\\
& \quad - \sum_{i=1}^N \Bigg( F^{\R}_i(\X^{\R}(t))
-4 \sum_{\substack{1 \leq j \leq N, \cr j \not=i}}
\frac{1}{X^{\R}_i(t)-X^{\R}_j(t)} \Bigg)
\frac{dt}{g_{\H^{\eta}_t}(z)-X^{\R}_i(t)}, 
\nonumber
\end{align}
\begin{align}
&d \Phi_{\O}(g_{\O^{\eta}_t}(z), \X^{\R_+}(t); q)
\label{eqn:dPhiO}
\\
& \qquad = - \sum_{i=1}^N 
\Bigg( \frac{1}{g_{\O^{\eta}_t}(z)-X^{\R_+}_i(t)}
- \frac{1}{g_{\O^{\eta}_t}(z)+X^{\R_+}_i(t)}
\Bigg) \sqrt{\kappa} d \widetilde{B}_i(t) 
\nonumber\\
& 
\qquad \quad
- \sum_{i=1}^N \Bigg[
F^{\R_+}_i(\X^{\R_+}(t)) 
\nonumber\\
& \qquad \qquad \quad
- \Bigg\{ 4 \sum_{\substack{1 \leq j \leq N, \cr j \not=i}}
\left( \frac{1}{X^{\R_+}_i(t)-X^{\R_+}_j(t)}
+ \frac{1}{X^{\R_+}_i(t)+X^{\R_+}_j(t)} \right) 
\nonumber\\
& \qquad \qquad \qquad \qquad \qquad \qquad \qquad 
\qquad \qquad 
+2 (1+2 \delta + q) \frac{1}{X^{\R_+}_i(t)}
\Bigg\}
\Bigg] 
\nonumber\\
& \qquad \qquad \qquad \qquad \qquad 
\times
\Bigg( \frac{1}{g_{\O^{\eta}_t}(z)-X^{\R_+}_i(t)}
- \frac{1}{g_{\O^{\eta}_t}(z)+X^{\R_+}_i(t)}
\Bigg) dt 
\nonumber\\
& 
\qquad \quad 
- 4 \delta \left( 1 - \frac{\kappa}{4} -q \right)
\frac{dt}{(g_{\O^{\eta}_{t}}(z))^2}
- \left(1 - \frac{\kappa}{4} \right) d \log g_{\O^{\eta}_t}'(z). 
\nonumber
\end{align}
\end{lem}
\noindent{\it Proof} \,
Apply It\^o's formula
and use the equalities such as
\[
\sum_{1 \leq i \not= j \leq N}
\frac{1}{(g-x_i)(g-x_j)}
=2 \sum_{1 \leq i \not= j \leq N}
\frac{1}{(g-x_i)(x_i-x_j)}.
\]
The proof is given by direct calculation.
\qed

If we assume that $(\X^{\R}(t))_{t \geq 0}$
is given by the $(8/\kappa)$-Dyson model $(\Y^{\R}(t))_{t \geq 0}$
satisfying (\ref{eqn:Dyson0}), the third term
in the RHS of (\ref{eqn:dPhiH}) vanishes.
Regarding (\ref{eqn:dPhiO}),  we first put
$q=1-\kappa/4$ to make the third term in the RHS become zero.
Then, if we assume that $\delta=\nu$ 
and $(\X^{\R_+}(t))_{t \geq 0}$
is given by the $(8/\kappa, \nu)$-Bru--Wishart process 
$(\Y^{\R_+}(t))_{t \geq 0}$
satisfying (\ref{eqn:Bru_WishartB}), the second term
in the RHS of (\ref{eqn:dPhiO}) vanishes.

Define 
\begin{align*}
\cM_{\H}(z, t)
&=- \Phi_{\H}(g_{\H^{\eta}_t}(z), \Y^{\R}(t))
- \left(1-\frac{\kappa}{4} \right)
\log g_{\H^{\eta}_t}'(z),
\nonumber\\
\cM_{\O}(z, t)
&=- \Phi_{\O}(g_{\O^{\eta}_t}(z), \Y^{\R_+}(t);  1-\kappa/4)
- \left(1-\frac{\kappa}{4} \right)
\log g_{\O^{\eta}_t}'(z).
\end{align*}

\begin{prop}
\label{thm:martingales}
Let $\kappa>0$,
$q=1-\kappa/4$, $\delta=\nu \geq 0$.
Then, for each point $z\in A$, 
$(\cM_D(z, t))_{0 \leq t \leq \tau}$,
$D=\H$ and $\O$, provide local martingales with increments
\begin{align*}
d \cM_{\H}(z, t)
&= \sum_{i=1}^N \frac{\sqrt{\kappa} dB_i(t)}
{g_{\H^{\eta}_t}(z)-Y^{\R}_i(t)},
\nonumber\\
d \cM_{\O}(z, t)
&= \sum_{i=1}^N 
\Bigg( \frac{1}{g_{\O^{\eta}_t}(z)-Y^{\R_+}_i(t)}
- \frac{1}{g_{\O^{\eta}_t}(z)+Y^{\R_+}_i(t)}
\Bigg) \sqrt{\kappa} d \widetilde{B}_i(t).
\end{align*}
\end{prop}

\subsection{GFF-valued Processes}
\label{sec:stationaryGFF}

Now, we consider the sum 
$H_{D^{\eta}_t}(\cdot)+ \sF[\cM_D(\cdot, t) ]$,
$0\leq t\leq \tau$,
where $\sF[\, \cdot \,]$
denotes a functional.
Comparing Lemma \ref{thm:dG} and 
Proposition \ref{thm:martingales} we observe that
\begin{equation}
d \Big\langle \Im \cM_D(z, \cdot), 
\Im \cM_D(w, \cdot) \Big\rangle_t
=- \frac{\kappa}{4} d G_{D^{\eta}_t}(z, w),
\label{eqn:dG}
\end{equation}
$z, w \in D^{\eta}_{\tau}$, $0 \leq t\leq \tau$, 
for $(D,S)=(\H, \R)$ and $(\O, \R_+)$.
Hence, we put 
$\sF[\, \cdot \,]=(2/\sqrt{\kappa}) \Im[\, \cdot \,]$,
and define the following GFF-valued processes
for $(D,S)=(\H, \R)$ and $(\O, \R_+)$, 
\begin{equation}
H_D(\cdot, t)
:= H_{D^{\eta}_t}(\cdot)+\frac{2}{\sqrt{\kappa}} 
\Im \cM_D(\cdot, t), \quad 0\leq t\leq \tau
\label{eqn:statGFF}
\end{equation}
with $\chi=\frac{2}{\sqrt{\kappa}}(1-\kappa/4)
=2/\sqrt{\kappa}-\sqrt{\kappa}/2$.
The second term of (\ref{eqn:statGFF}) contains 
an imaginary part of the complex-valued
logarithmic potential
$-\Phi_D(g_{D^{\eta}_t}(z), \Y^S(t))$, $t \geq 0$.
This is the unique harmonic function satisfying
the boundary condition
\[
\frac{2}{\sqrt{\kappa}} \Im \cM_D(x, t)
=\begin{cases}
\displaystyle{
-\frac{2 \pi}{\sqrt{\kappa}} N},
&\mbox{if $x < Y^S_1(t)$},
\cr
\cr
\displaystyle{
-\frac{2 \pi}{\sqrt{\kappa}}(N-i)},
&\mbox{if $x \in (Y^S_i(t), Y^S_{i+1}(t)), 1 \leq i \leq N$},
\end{cases}
\]
with the convention $Y^{S}_{N+1}(t) \equiv + \infty$.
That is, it has discontinuity at $Y^S_i(t)$ by 
$2 \pi/\sqrt{\kappa}$ along $S$, $1 \leq i \leq N$, $t \geq 0$.
We will think that the GFF $H_D(\cdot, t)$ has 
the same boundary condition as 
$(2/\sqrt{\kappa}) \Im \cM_D(\cdot, t)$, 
$0\leq t\leq \tau$.
For further arguments concerning 
the second term of (\ref{eqn:statGFF}),
see Section V.C in \cite{KK20}.

\begin{thm}
\label{thm:main1}
Let $\kappa>0$, $q=1-\kappa/4, \delta=\nu \geq 0$.
Assume that $(D, S)=(\H, \R)$ or $(\O, \R_+)$, 
and $(\Y^S(t))_{t \geq 0}$ is the $(8/\kappa)$-Dyson model
if $S=\R$ and the $(8/\kappa, \nu)$-Bru--Wishart process
if $S=\R_+$,
starting from a configuration in $\W_N(S)$.
Then, for each $f \in \cC_{\rm c}^{\infty}(D)$ 
such that $\supp (f)\subset A$,
we have the following equality:
\[
(H_D(\cdot, 0), f)
\law= (H_D(\cdot, \tau), f )
\quad \mbox{under $\P$}.
\]
\end{thm}
\noindent{\it Proof} \,
From (\ref{eqn:dG}) we have
\begin{equation}
	d \left\langle 
  \Big( \frac{2}{\sqrt{\kappa}} 
	\Im \cM_D (\cdot,\cdot), f \Big) \right\rangle_{t}
	=-dE_{t}(f),\ \ 0\le t\le \tau,
\label{eqn:dE}
\end{equation}
where
\begin{equation}
E_{t}(f):=\int_{A\times A}f(z)G_{D^{\eta}_{t}}(z,w)f(w)d\mu (z)d\mu (w) 
\label{eqn:D_energy}
\end{equation}
is called the \textit{Dirichlet energy} 
of $f$ in $D^{\eta}_{t}$. 
We have $\Var [ (H_{D^{\eta}_{t}},f)]=E_{t}(f)$
due to the conformal invariance of a zero-boundary GFF.
Introducing a parameter $\theta \in \R$, we see 
\[
\E \left[e^{\sqrt{-1} \theta (H_{D}(\cdot,\tau),f)} \right]
=\E \left[\E 
\left[e^{\sqrt{-1} \theta(H_{D^{\eta}_{\tau}},f)}
\Big|\cF_{\tau}\right] 
e^{\sqrt{-1} \theta\frac{2}{\sqrt{\kappa}}
(\Im \cM_{D}(\cdot,\tau),f)} \right],
\]
where $\E$ is the expectation value 
with respect to the probability measure $\P$,
since $(\Im \cM_{D}(\cdot,\tau),f)$ is 
$\cF_{\tau}$-measurable. 
By definition of a zero-boundary GFF
and the Dirichlet energy (\ref{eqn:D_energy}), 
we obtain
\[
\E\left[e^{\sqrt{-1} \theta(H_{D^{\eta}_{\tau}},f)}
\Big|\cF_{\tau}\right]=e^{-\frac{\theta^{2}}{2}
E_{\tau}(f)}. 
\]
Hence, by Proposition \ref{thm:martingales}
and (\ref{eqn:dE}) we have
\begin{align*}
\E \left[e^{\sqrt{-1} \theta (H_{D}(\cdot,\tau),f)}\right]
&=\E \left[e^{\sqrt{-1} \theta\frac{2}{\sqrt{\kappa}}
(\Im \cM_{D}(\cdot,\tau),f)-\frac{\theta^{2}}{2}
E_{\tau}(f)}\right] 
\nonumber\\
&=\E \left[ e^{\sqrt{-1} \theta\frac{2}{\sqrt{\kappa}}
(\Im \cM_{D}(\cdot,0),f)-\frac{\theta^{2}}{2}E_{0}(f)} \right]
\nonumber\\
&=\E \left[e^{\sqrt{-1} \theta (H_{D}(\cdot,0),f)}\right].
\end{align*}
This implies the desired coincidence 
under the probability law $\P$. 
\qed

\vskip 0.3cm
\noindent{\bf Remark 5.1} 
There are two local formulations of multiple SLE; 
one of them is based on commutation relation of Loewner chains 
each of which generates a single random curve \cite{Dub07}
and the other one is a single Loewner chain 
driven by multiple driving processes \cite{BBK05}. 
Though it is expected that these two formulations are equivalent, 
the equivalence has not been proved and even its precise 
statement is not obvious. 
In the work by Miller and Sheffield \cite{MS16a}, 
they studied coupling between a GFF and a variant of SLE 
called SLE($\kappa, \Urho$), which reduces to 
a member of commuting Loewner chains at a specific 
setting of parameters. Hence, it can be said that they 
also considered coupling between a GFF and 
a multiple SLE in the former sense. 
This does not, however, imply that a multiple SLE 
in the latter sense can be coupled with the same GFF, 
which is exactly the result we presented in the current article.

\vskip 1cm
\noindent{\bf Acknowledgements} \,
The present authors would like 
to thank Kalle Kyt\"{o}l\"{a}
for useful comments on multiple SLEs.
MK was supported by
the Grant-in-Aid for Scientific Research (C) (No.19K03674),
(B) (No.18H01124), 
and (S) (No. 16H06338) 
of Japan Society for the Promotion of Science (JSPS).
SK was supported by the Grant-in-Aid
for JSPS Fellows (No.19J01279). 


\end{document}